\documentclass[12pt]{article}
\usepackage{amssymb,amsmath}
\usepackage{cases}
\usepackage{amsfonts}
\usepackage{cite,color,xcolor}
\usepackage[left=2.0cm,right=2.0cm,top=2.0cm,bottom=2.0cm]{geometry}
\usepackage[colorlinks,citecolor=blue,urlcolor=blue]{hyperref}
\usepackage[utf8]{inputenc}

\newtheorem{theorem}{Theorem}[section]

\newtheorem{lemma}{Lemma}[section]
\newtheorem{proposition}{Proposition}[section]

\newtheorem{definition}{Definition}[section]
\newtheorem{remark}{Remark}[section]

\usepackage{txfonts}

\newcommand{\bal}{\begin{align}}
\newcommand{\bbal}{\begin{align*}}
\newcommand{\beq}{\begin{equation}}
\newcommand{\eeq}{\end{equation}}
\newcommand{\bca}{\begin{cases}}
\newcommand{\eca}{\end{cases}}
\def\div{\mathord{{\rm div}}}
\newcommand{\pa}{\partial}
\newcommand{\fr}{\frac}
\newcommand{\na}{\nabla}
\newcommand{\De}{\Delta}

\newcommand{\cd}{\cdot}
\newcommand{\ep}{\varepsilon}
\newcommand{\dd}{\mathrm{d}}

\newcommand{\R}{\mathbb{R}}

\newcommand{\les}{\lesssim}

\newcommand{\bi}{\Big}
\newcommand{\g}{\big}
\begin{document}
\title{Ill-posedness for the Euler equations in Besov spaces}

\author{Jinlu Li$^{1}$, Yanghai Yu$^{2,}$\footnote{E-mail: lijinlu@gnnu.edu.cn; yuyanghai214@sina.com(Corresponding author); mathzwp2010@163.com} and Weipeng Zhu$^{3}$\\
\small $^1$ School of Mathematics and Computer Sciences, Gannan Normal University, Ganzhou 341000, China\\
\small $^2$ School of Mathematics and Statistics, Anhui Normal University, Wuhu 241002, China\\
\small $^3$ School of Mathematics and Big Data, Foshan University, Foshan, Guangdong 528000, China}

\date{\today}

\maketitle\noindent{\hrulefill}

{\bf Abstract:} In the paper, we consider the Cauchy problem to the Euler equations in $\R^d$ with $d\geq2$. We construct an initial data $u_0\in B^\sigma_{p,\infty}$ showing that the corresponding solution map of the Euler equations starting from $u_0$ is
discontinuous at $t = 0$ in the metric of $B^\sigma_{p,\infty}$, which implies the ill-posedness for this equation in $B^\sigma_{p,\infty}$. We generalize the periodic result of Cheskidov and Shvydkoy \cite{Cheskidov} (Proc. Amer. Math. Soc. 138 (2010), 1059--1067).


\vskip0mm\noindent{\hrulefill}

\section{Introduction}

In this article, we consider the Euler equations governing the motion of an incompressible
fluid in $\R^d$ with $d\geq2$
\begin{align}\label{e}
\begin{cases}
\pa_t u+u\cdot \nabla u+\nabla P=0, &\quad (t,x)\in \R^+\times\R^d,\\
\mathrm{div\,} u=0,&\quad (t,x)\in \R^+\times\R^d,\\
u(0,x)=u_0(x), &\quad x\in \R^d,
\end{cases}
\end{align}
where the vector field $u(t,x):[0,\infty)\times {\mathbb R}^d\to {\mathbb R}^d$ stands for the velocity of the fluid, the quantity $P(t,x):[0,\infty)\times {\mathbb R}^d\to {\mathbb R}$ denotes
the scalar pressure, and $\mathrm{div\,} u=0$ means that the fluid is incompressible. The mathematical study of the Euler equations of ideal hydrodynamics has a long and distinguished history. We do not detail the literature since it is huge and refer the readers to see
the monographs of Bahouri--Chemin--Danchin \cite{B.C.D} and Majda--Bertozzi \cite{Majda}.

We say that the Cauchy problem \eqref{e} is Hadamard (locally) well-posed in a Banach space $X$ if for any data $u_0\in X$ there exists (at least for a short time) $T>0$ and a unique solution in the space $\mathcal{C}([0,T),X)$ which depends continuously on the data. In particular, we say that the solution map is continuous if for any $u_0\in X$, there exists a neighborhood $B \subset X$ of $u_0$ such that
for every $u \in B$ the map $u \mapsto U$ from $B$ to $\mathcal{C}([0, T]; X)$ is continuous, where $U$
denotes the solution to \eqref{e} with initial data $u_0$. Our goal of this paper is to consider the ill-posedness of \eqref{e} in $B^\sigma_{p,\infty}(\R^d)$ and generalize the periodic result of Cheskidov and Shvydkoy \cite{Cheskidov}.
Our main result of this paper reads as follows:
\begin{theorem}\label{th1}
Let $d\geq2$. Assume that  $$\sigma>1+\fr{d}{p}\quad \text{with}\quad 1\leq p\leq \infty.$$ There exits $u_0\in B^\sigma_{p,\infty}(\R^d)$ and a positive constant $\ep_0$ such that the data-to-solution map $u_0\mapsto \mathbf{S}_{t}(u_0)$ of the Cauchy problem \eqref{e} satisfies\bbal
\limsup_{t\to0^+}\|\mathbf{S}_{t}(u_0)-u_0\|_{B^\sigma_{p,\infty}}\geq \ep_0.
\end{align*}
\end{theorem}
\begin{remark}\label{r2}
Theorem \ref{th1} demonstrates the ill-posedness of the Euler equations in $B^\sigma_{p,\infty}$. More precisely, there exists $u_0\in B^\sigma_{p,\infty}$ such that the corresponding solution to the Euler equations that starts from $u_0$ does not converge back to $u_0$ in the metric of $B^\sigma_{p,\infty}$ as time goes to zero.
\end{remark}
\begin{remark}\label{r1}
The method we used in proving Theorem \ref{th1} is completely different from that \cite{Cheskidov}.
\end{remark}
\begin{remark}\label{r3}
Due to the incompressible condition $\eqref{e}_2$, the pressure can be eliminated from \eqref{e}. In fact, applying the Leray operator $\mathcal{P}$ to $\eqref{e}_1$, then we have
\begin{align*}
\pa_t u=-\mathcal{P}(u\cdot \nabla u)=:\mathbf{P}(u).
\end{align*}
Thus, we transform \eqref{e} equivalently into the nonlocal form
\begin{align}\label{ch}
\begin{cases}
\pa_t u=\mathbf{P}(u), \\
\mathrm{div\,} u=0,\\
u(0,x)=u_0(x).
\end{cases}
\end{align}
\end{remark}

Concerning the well-posedness of the Euler equations we refer to see \cite{B.C.D,guo,Majda}. Next, we recall some results about the properties of the data-to-solution map. The first results of this type was proved by Kato \cite{Kato} who showed that the
solution operator for the (inviscid) Burgers equation is not H\"{o}lder continuous in the
$H^s(\mathbb{T})$-norm $(s > 3/2)$ for any H\"{o}lder exponent. Kato--Ponce \cite{Kato87} obtained the continuity results for the solution map of the Euler equations in $W^s_p$ (see also \cite{Kato84}).
Pak--Park \cite{Pak} established existence and uniqueness
of solutions of the Euler equations in $B^{1}_{\infty,1}$ and showed that the solution map is in fact Lipschitz
continuous when viewed as a map between $B^{0}_{\infty,1}$ and $\mathcal{C}([0,\infty);B^{0}_{\infty,1})$. Later, Cheskidov--Shvydkoy \cite{Cheskidov} proved that the solution of the Euler equations cannot be continuous as a function
of the time variable at $t = 0$ in the spaces $B^s_{r,\infty}(\mathbb{T}^d)$ where $s > 0$ if $2 < r \leq \infty$ and $s>d(2/r-1)$
if $1 \leq r \leq 2$. Furthermore, Bourgain--Li \cite{Bourgain1,Bourgain2} proved the strong local ill-posedness of the Euler equations in borderline Besov spaces $B^{d/p+1}_{p,r}$ with $(p,r)\in[1,\infty)\times(1,\infty]$ when $d=2,3$. Subsequently, Misio{\l}ek--Yoneda \cite{MY} studied the borderline cases and showed that the 2D Euler equations are not locally well-posed in the sense of Hardamard in the $C^1$ space and in the Besov space $B^1_{\infty,1}$.
Recently, Misio{\l}ek--Yoneda \cite{MY2} showed that the solution map for the Euler equations is not even continuous in the space of H\"{o}lder continuous functions and thus not locally Hadamard well-posed in $B^{1+\sigma}_{\infty,\infty}$ with any $\sigma\in(0,1)$. Concerning the non-uniform continuity of the data-to-solution map, we refer the readers to see \cite{H1,H2,H3,Pastrana,Li}.

\section{Littlewood-Paley Decomposition and Besov Spaces}
We will use the following notations throughout this paper. The notation $A\les B$ (resp., $A \gtrsim B$) means that there exists a harmless positive constant $c$ such that $A \leq cB$ (resp., $A \geq cB$). For $X$ a Banach space and $I\subset\R$, we denote by $\mathcal{C}(I;X)$ the set of continuous functions on $I$ with values in $X$.
Let us recall that for all $u\in \mathcal{S}'$, the Fourier transform $\mathcal{F}u$, also denoted by $\widehat{u}$, is defined by
$$
\mathcal{F}u(\xi)=\widehat{u}(\xi)=\int_{\R^d}e^{-ix\cd \xi}u(x)\dd x \quad\text{for any}\; \xi\in\R^d.
$$
Due to the Holdge decomposition, we know that any vector field $f= (f_1,...,f_d)$ with components in $S_h'(\mathbb{R}^d)$ may be decomposed into one potential part $\mathcal{Q}f$ and one divergence-free part $\mathcal{P}f$, where the projectors $\mathcal{P}$ and $\mathcal{Q}$ are defined by
\begin{eqnarray}\label{Equ2.1}
\mathcal{P}= \mathbb{I}+(-\Delta)^{-1}\nabla {\rm{div}} \quad\mbox{and}\quad \mathcal{Q}= -(-\Delta)^{-1}\nabla {\rm{div}}.
\end{eqnarray}

Next, we will recall some facts about the Littlewood-Paley decomposition, the nonhomogeneous Besov spaces and their some useful properties.
\begin{proposition}[Littlewood-Paley decomposition, See \cite{B.C.D}] Let $\mathcal{B}:=\{\xi\in\mathbb{R}^d:|\xi|\leq \frac 4 3\}$ and $\mathcal{C}:=\{\xi\in\mathbb{R}^d:\frac 3 4\leq|\xi|\leq \frac 8 3\}.$
There exist two radial functions $\chi\in C_c^{\infty}(\mathcal{B})$ and $\varphi\in C_c^{\infty}(\mathcal{C})$ both taking values in $[0,1]$ such that
\begin{align*}
&\chi(\xi)+\sum_{j\geq0}\varphi(2^{-j}\xi)=1 \quad \forall \;  \xi\in \R^d,\\
&\frac{1}{2} \leq \chi^{2}(\xi)+\sum_{j \geq 0} \varphi^{2}(2^{-j} \xi) \leq 1\quad \forall \;  \xi\in \R^d.
\end{align*}
\end{proposition}
Moreover, we can deduce that $\varphi\equiv 1$ for $\frac43\leq |\xi|\leq \frac32$. This basic fact will be used in the sequel.

For every $u\in \mathcal{S'}(\mathbb{R}^d)$, the inhomogeneous dyadic blocks ${\Delta}_j$ are defined as follows
\begin{numcases}{\Delta_ju=}
0, & if $j\leq-2$;\nonumber\\
\chi(D)u=\mathcal{F}^{-1}(\chi \mathcal{F}u), & if $j=-1$;\nonumber\\
\varphi(2^{-j}D)u=\mathcal{F}^{-1}\g(\varphi(2^{-j}\cdot)\mathcal{F}u\g), & if $j\geq0$.\nonumber
\end{numcases}
In the inhomogeneous case, the following Littlewood-Paley decomposition makes sense
$$
u=\sum_{j\geq-1}{\Delta}_ju\quad \text{for any}\;u\in \mathcal{S'}(\mathbb{R}^d).
$$

We turn to the definition of the Besov Spaces and norms which will come into play in our paper.
\begin{definition}[Besov Space \cite{B.C.D}]
Let $s\in\mathbb{R}$ and $(p,r)\in[1, \infty]^2$. The nonhomogeneous Besov space $B^s_{p,r}(\R^d)$ consists of all tempered distribution $f$ such that
\begin{align*}
B^{s}_{p,r}(\R^d):=\Big\{f\in \mathcal{S}'(\mathbb{R}):\;\|f\|_{B^{s}_{p,r}(\R^d)}<\infty\Big\},
\end{align*}
where
\begin{numcases}{\|f\|_{B^{s}_{p,r}(\R^d)}:=}
\left(\sum_{j\geq-1}2^{sjr}\|\Delta_jf\|^r_{L^p(\R^d)}\right)^{1/r}, &if $1\leq r<\infty$,\nonumber\\
\sup_{j\geq-1}2^{sj}\|\Delta_jf\|_{L^p(\R^d)}, &if $r=\infty$.\nonumber
\end{numcases}
\end{definition}
\begin{remark}\label{re3}
It should be emphasized that the following facts will be often used implicity:
\begin{itemize}
  \item $B^s_{p,q}(\R)\hookrightarrow B^t_{p,r}(\R)\quad\text{for}\;s>t\quad\text{or}\quad s=t,1\leq q\leq r\leq\infty.$
  \item For $1\leq p\leq\infty$ and $s>\frac{d}{p}+1$, $B^s_{p,\infty}(\R^d)$ is a Banach algebra.
\end{itemize}
\end{remark}
\section{Proof of Theorem \ref{th1}}\label{sec3}
In this section, we prove Theorem \ref{th1} only for the case $d=2$ since the case $d\geq3$ requires very little changes (for more details see \cite{Li}).
\subsection{Construction of Initial Data}
Firstly, we need to introduce smooth, radial cut-off functions to localize the frequency region.
Let $\widehat{\phi}\in \mathcal{C}^\infty_0(\mathbb{R})$ be an even, real-valued and non-negative function on $\R$ and satisfy
\begin{numcases}{\widehat{\phi}(\xi)=}
1, &if $|\xi|\leq \frac{1}{16}$,\nonumber\\
0, &if $|\xi|\geq \frac{1}{4}$.\nonumber
\end{numcases}

We set the initial data $u_0(x)$ as follows:
\bal\label{cz}
u_0(x):=\sum\limits^\infty_{j=0}f_{kj}(x),\quad x\in\R^2,\;k\in\mathbb{N}^+,
\end{align}
where
\bbal
f_{m}(x):=2^{-m(\sigma+1)}
\begin{pmatrix}
-\pa_2\\
\pa_1
\end{pmatrix}
g_{m} \quad\text{with}\quad g_{m}(x):=\phi(x_1)\cos \bi(\frac{17}{12}2^{m}x_1\bi)\phi(x_2),\quad m\geq0.
\end{align*}
\subsection{A Key Lemma}

We establish the following crucial Lemmas which involves the properties of the initial data.
\begin{lemma}\label{ley2}
Define the initial data $u_0(x)$ by \eqref{cz} above. Then there exist some sufficiently large $k$ and $n$ such that
\begin{description}
  \item[(1)] $\div u_0=0$;
  \item[(2)] $\|u_0\|_{B^\sigma_{p,\infty}}\leq C$;
  \item[(3)] $\|\De_{kn}(u_0\cd\na u_0)\|_{L^p}\geq c2^{kn}2^{-kn\sigma}$,
\end{description}
where $C$ and $c$ are uniform constants.
\end{lemma}
{\bf Proof.}\quad A direct calculation yields {\bf(1)}.

Easy computations give that
\bbal
\widehat{f}_{m}(\xi)=2^{-{m}(\sigma+1)-1}
\begin{pmatrix}
-i\xi_2\\
i\xi_1
\end{pmatrix}
\bi[\widehat{\phi}\bi(\xi_1-\frac{17}{12}2^{m}\bi)+\widehat{\phi}\bi(\xi_1+\frac{17}{12}2^{m}\bi)\bi]
\widehat{\phi}(\xi_2),
\end{align*}
which implies
\bbal
\mathrm{supp} \ \widehat{f}_{m}&\subset \Big\{\xi\in\R^2: \ \frac{17}{12}2^{m}-\fr12\leq |\xi|\leq \frac{17}{12}2^{m}+\fr12\Big\}\\
&\subset \Big\{\xi\in\R^2: \ \frac{33}{24}2^{m}\leq |\xi|\leq \frac{35}{24}2^{m}\Big\}\quad\text{for}\quad m\geq4.
\end{align*}
Notice that
\bbal
\varphi(2^{-j}\xi)\equiv 1\quad \text{in}\;  \Big\{\xi\in\R^2: \ \frac{4}{3}2^{j}\leq |\xi|\leq \frac{3}{2}2^{j}\Big\},
\end{align*}
and
\bbal
\mathcal{F}\big(\Delta_j(f_{m})\big)=\varphi(2^{-j}\cdot)\widehat{f}_{m},
\end{align*}
then we deduce for $m\geq4$
\begin{align}\label{y}
\Delta_j(f_{m})=\begin{cases}
f_m, &if\;\; j=m,\\
0, &if\;\; j\neq m.
\end{cases}
\end{align}
We should mention that \eqref{y} also holds if $\emph{cos}$ were replaced by $\emph{sin}$ in the definition of $g_m$.

Thus, we deduce from \eqref{y} that
\bbal
\|u_0\|_{B^\sigma_{p,\infty}(\R^2)}&=\sup_{\ell\geq-1}2^{\ell\sigma}\|\Delta_\ell u_0\|_{L^{p}(\R^2)}
\leq C.
\end{align*}
This gives {\bf(2)}.

To prove {\bf(3)}, we write
\bal\label{yy}
u_0\cd\na u_0=\sum^\infty_{j=0}f_{kj}\cd\na f_{kj}+\sum^\infty_{j=1}\sum^{j-1}_{i=0}(f_{kj}\cd\na f_{ki}+f_{ki}\cd\na f_{kj}).
\end{align}
On the one hand, we have for $\ell=1,2$
\bbal
\g(f_{kj}\cd\na f_{kj}\g)_\ell
&=2^{-2kj(\sigma+1)}\Big[\Phi_1(x_1, x_2)+\Phi_2(x_1, x_2)\cos \bi(\frac{17}{12}2^{kj+1}x_1\bi)+\Phi_3(x_1, x_2)\sin \bi(\frac{17}{12}2^{kj+1}x_1\bi)\Big],
\end{align*}
where
$$\mathrm{supp}\ \widehat{\Phi_1},\;\mathrm{supp}\ \widehat{\Phi_2},\;\mathrm{supp}\ \widehat{\Phi_3}\subset \big\{\xi\in\R^2: \ |\xi|\leq 2\big\}.$$
Then we have
\bbal
\Delta_{kn}\left((f_{kj}\cd\na f_{kj})_\ell\right)\equiv0\quad\text{for}\quad j\geq0,
\end{align*}
which implies that
\bal\label{yy1}
\Delta_{kn}\left(\sum^\infty_{j=0}f_{kj}\cd\na f_{kj}\right)\equiv0.
\end{align}
On the other hand, notice that
\bbal
\mathrm{supp}\ \mathcal{F}\left(f_{kj}\cd \na f_{ki}\right)
&\subset \left\{\xi\in\R^2: \ \frac{17}{12}2^{kj}-\frac{17}{12}2^{ki}-\fr12\leq |\xi|\leq \frac{17}{12}2^{kj}+\frac{17}{12}2^{ki}+\fr12\right\}\\
&\subset \left\{\xi\in\R^2: \ \frac{33}{24}2^{kj}\leq |\xi|\leq \frac{35}{24}2^{kj}\right\},
\end{align*}
and \eqref{y}, we have
\bal\label{yy2}
\Delta_{kn}\left(\sum^\infty_{j=1}\sum^{j-1}_{i=0}(f_{kj}\cd\na f_{ki}+f_{ki}\cd\na f_{kj})\right)=\sum^{n-1}_{i=0}(f_{kn}\cd\na f_{ki}+f_{ki}\cd\na f_{kn}).
\end{align}
Thus, we deduce from \eqref{yy1} and \eqref{yy2} that
\bbal
\Delta_{kn}(u_0\cd\na u_0)&=\sum^{n-1}_{i=0}f_{kn}\cd\na f_{ki}+\sum^{n-1}_{i=1}f_{ki}\cd\na f_{kn}+f_{0}\cd\na f_{kn}\\
&=:\mathbf{I}_1+\mathbf{I}_2+\mathbf{I}_3.
\end{align*}
The first two terms can be estimated as follows:
\bal
\|\mathbf{I}_1\|_{L^p}&\leq\sum^{n-1}_{i=0}\|f_{kn}\cd\na f_{ki}\|_{L^p}\nonumber\\
&\leq\|f_{kn}\|_{L^p}\sum^{n-1}_{i=0}\|\na f_{ki}\|_{L^\infty}\nonumber\\
&\leq  C2^{-kn\sigma}\sum^{n-1}_{i=0}2^{-ki(\sigma-1)}\nonumber\\
&\leq C2^{-kn\sigma}\label{l1}
\end{align}
and
\bal
\|\mathbf{I}_2\|_{L^p}&\leq\sum^{n-1}_{i=1}\|f_{ki}\cd\na f_{kn}\|_{L^p}\nonumber\\
&\leq\|\nabla f_{kn}\|_{L^p}\sum^{n-1}_{i=1}\|f_{ki}\|_{L^\infty}\nonumber\\
&\leq C2^{kn}2^{-kn\sigma}\sum^{n-1}_{i=1}2^{-k\sigma i}\nonumber\\
&\leq C2^{kn}2^{-kn\sigma}2^{-k\sigma}.\label{l2}
\end{align}
For the last term $f_{0}\cd\na f_{kn}$, we compute the two components:
\bbal
\big(f_{0}\cd\na f_{kn})_1&=-2^{-kn(\sigma+1)}\big[(f_0)_1\pa_1\pa_2g_{kn}+(f_0)_2\pa_2^2g_{kn}\big]\\
&=2^{-kn(\sigma+1)}(\pa_2g_0\pa_1\pa_2g_{kn}-\pa_1g_0\pa_2^2g_{kn})\\
&=2^{-kn(\sigma+1)}(\pa_2g_0h_2-\pa_1g_0h_1)=:\mathbf{I}_{3,1}
\end{align*}
and
\bbal
\big(f_{0}\cd\na f_{kn})_2&=2^{-kn(\sigma+1)}\big[(f_0)_1\pa^2_1g_{kn}+(f_0)_2\pa_1\pa_2g_{kn}\big]
\\&=2^{-kn(\sigma+1)}(\pa_1g_0h_2-\pa_2g_0h_3)+\Big(\frac{17}{12}\Big)^22^{-kn(\sigma-1)}h_4\\&=:\mathbf{I}_{3,2}+\mathbf{I}_{3,3},
\end{align*}
where
\bbal
h_1&:=\phi(x_1)\cos \bi(\frac{17}{12}2^{kn}x_1\bi)\phi''(x_2),\\
h_2&:=\phi'(x_1)\cos \bi(\frac{17}{12}2^{kn}x_1\bi)\phi'(x_2)-\frac{17}{12}2^{kn}\phi(x_1)\sin \bi(\frac{17}{12}2^{kn}x_1\bi)\phi'(x_2),\\
h_3&:=\phi''(x_1)\cos \bi(\frac{17}{12}2^{kn}x_1\bi)\phi(x_2)-\frac{17}{6}2^{kn}\phi'(x_1)\sin \bi(\frac{17}{12}2^{kn}x_1\bi)\phi(x_2),\\
h_4&:=\phi^2(x_1)\cos \bi(\frac{17}{12}x_1\bi)\cos \bi(\frac{17}{12}2^{kn}x_1\bi)\phi(x_2)\phi'(x_2).
\end{align*}
By simple computations, one has
\bbal
\|\mathbf{I}_{3,1}\|_{L^p}+\|\mathbf{I}_{3,2}\|_{L^p}&\leq 2^{-kn(\sigma+1)}\|\pa_1g_0,\pa_2g_0\|_{L^p}\|h_1,h_2,h_3\|_{L^\infty}\leq C2^{-kn\sigma}
\end{align*}
and
\bbal
&\|\mathbf{I}_{3,3}\|_{L^p}=\Big(\frac{17}{12}\Big)^22^{-kn(\sigma-1)}\|h_4\|_{L^p}\geq c2^{-kn(\sigma-1)},
\end{align*}
where we have used the simple fact for sufficiently large $n$
$$\|h_4\|_{L^p}=\bi\|\phi^2(x_1)\cos \bi(\frac{17}{12}x_1\bi)\cos \bi(\frac{17}{12}2^nx_1\bi)\bi\|_{L^p(\R)}\bi\|\phi(x_2)\phi'(x_2)\bi\|_{L^p(\R)}\geq c_0.
$$
This gives
\bal\label{l3}
\|\mathbf{I}_{3}\|_{L^p}\geq c_02^{-kn(\sigma-1)}-C2^{-kn\sigma}.
\end{align}
Combining \eqref{l1}--\eqref{l3} yields
\bbal
\|\Delta_{kn}(u_0\cd\na u_0)\|_{L^p}\geq c_0\cdot 2^{kn}2^{-kn\sigma}-C2^{kn}2^{-kn\sigma}2^{-k\sigma}-C2^{-kn\sigma}.
\end{align*}
 We take sufficiently large $k$ and $n$ and obtain the desired {\bf(3)}. Thus, we complete the proof of Lemma \ref{ley2}.
\subsection{Error Estimates}\label{sec3.2}

\begin{proposition}\label{pro3.2}
Let $u_0\in B^{\sigma}_{p,\infty}$. Assume that $u\in L^\infty_TB^{\sigma}_{p,\infty}$ be the solution of the Cauchy problem \eqref{ch}, we have
\bbal
\|\mathbf{w}(t,u_0)\|_{\dot{B}^{\sigma-2}_{p,\infty}}\leq Ct^2\|u_0\|^3_{B^\sigma_{p,\infty}},
\end{align*}
here and in what follows we denote
$$\mathbf{w}(t,u_0):=\mathbf{S}_{t}(u_0)-u_0-t\mathbf{P}(u_0).$$
\end{proposition}
{\bf Proof.}\quad For simplicity, we denote $u(t):=\mathbf{S}_t(u_0)$ here and in what follows. Due to the fact $B^\sigma_{p,\infty}\hookrightarrow\rm Lip$, we know that there exists a positive time $T=T(\|u_0\|_{B^\sigma_{p,\infty}})$ such that
\bal\label{s}
\|u(t)\|_{L^\infty_TB^\sigma_{p,\infty}}\leq C\|u_0\|_{B^\sigma_{p,\infty}}\leq C.
\end{align}
Moreover, for $\gamma\geq \sigma-1$, we have
\bal\label{s1}
\|u(t)\|_{L^\infty_TB^\gamma_{p,\infty}}\leq C\|u_0\|_{B^\gamma_{p,\infty}}.
\end{align}
By the Mean Value Theorem, we obtain from \eqref{ch} that
\bal\label{j}
\|u(t)-u_0\|_{B^{\sigma-1}_{p,\infty}}
&\leq \int^t_0\|\pa_\tau u\|_{B^{\sigma-1}_{p,\infty}} \dd\tau
\nonumber\\
&\leq \int^t_0\|\mathbf{P}(u)\|_{B^{\sigma-1}_{p,\infty}} \dd\tau
\nonumber\\
&\les \int^t_0\|\mathbf{P}(u)\|_{\dot{B}^{\sigma-1}_{p,\infty}\cap\dot{B}^0_{p,1}} \dd\tau
\nonumber\\
&\les \int^t_0\|u\otimes u\|_{\dot{B}^{\sigma}_{p,\infty}\cap\dot{B}^1_{p,1}} \dd\tau
\nonumber\\
&\les \int^t_0\|u\otimes u\|_{B^{\sigma}_{p,\infty}} \dd\tau\nonumber\\
&\les t\|u_0\|_{B^{\sigma}_{p,\infty}}^2,
\end{align}
where we have used the fact that $B^\sigma_{p,\infty}$ is a Banach algebra.

By the Mean Value Theorem and \eqref{ch}, we obtain
\bbal
\|\mathbf{w}(t,u_0)\|_{\dot{B}^{\sigma-2}_{p,\infty}}
&\leq \int^t_0\|\pa_\tau u-\mathbf{P}(u_0)\|_{\dot{B}^{\sigma-2}_{p,\infty}} \dd\tau
\\&\leq \int^t_0\|\mathbf{P}(u)-\mathbf{P}(u_0)\|_{\dot{B}^{\sigma-2}_{p,\infty}} \dd\tau
\\
&\les \int^t_0\|u\otimes u-u_0\otimes u_0\|_{\dot{B}^{\sigma-1}_{p,\infty}} \dd\tau
\\
&\les \int^t_0\|u\otimes u-u_0\otimes u_0\|_{{B}^{\sigma-1}_{p,\infty}} \dd\tau
\\&\les\int^t_0\|u(\tau)-u_0\|_{B^{\sigma-1}_{p,\infty}} \|u_0,u(\tau)\|_{B^{\sigma}_{p,\infty}} \dd\tau
\\&\les t^2\|u_0\|^3_{B^{\sigma}_{p,\infty}},
\end{align*}
where we have used \eqref{j} in the last step.

Thus, we complete the proof of Proposition \ref{pro3.2}.

Now we present the proof of Theorem \ref{th1}.\\
{\bf Proof of Theorem \ref{th1}.}\quad
Using Proposition \ref{pro3.2} and Lemma \ref{ley2}, we have
\bbal
\|\mathbf{S}_{t}(u_0)-u_0\|_{B^\sigma_{p,\infty}}
&\geq2^{{kn\sigma}}\big\|\De_{kn}\big(\mathbf{S}_{t}(u_0)-u_0\big)\big\|_{L^p}\\
&=2^{{kn\sigma}}\big\|\De_{kn}\big(t\mathbf{P}(u_0)+\mathbf{w}(t,u_0)\big)\big\|_{L^p}\\
&\geq t2^{{kn\sigma}}\big\|\De_{kn}\big(\mathcal{P}(u_0\cd\na u_0)\big)\big\|_{L^p}-2^{{kn\sigma}}\big\|\De_{kn}\big(\mathbf{w}(t,u_0)\big)\big\|_{L^p}\\
&\geq t2^{{kn\sigma}}\|\De_{kn}\big(u_0\cd\na u_0\big)\|_{L^p}-t2^{{kn\sigma}}\|\De_{kn}\mathcal{Q}(u_0\cd\na u_0)\|_{L^p}\\
&~~~~-2^{{2kn}}2^{{kn(\sigma-2)}}\big\|\De_{kn}\big(\mathbf{w}(t,u_0)\big)\big\|_{L^p}\\
&\geq t2^{{kn\sigma}}\|\De_{kn}\big(u_0\cd\na u_0\big)\|_{L^p}-
Ct\|\mathcal{Q}(u_0\cd\na u_0)\|_{B^{\sigma}_{p,\infty}}\\
&~~~~-C2^{2{kn}}\|\mathbf{w}(t,u_0)\|_{\dot{B}^{\sigma-2}_{p,\infty}}\\
&\geq ct2^{{kn}}-Ct-C2^{2{kn}}t^2,
\end{align*}
where we have used the fact:
\bbal
\|\mathcal{Q}(u_0\cd\na u_0)\|_{B^{\sigma}_{p,\infty}}= \|(-\Delta)^{-1}\nabla(\pa_iu^j_0\pa_j u_0^i)\|_{B^{\sigma}_{p,\infty}}\les \|\pa_iu^j_0\pa_j u_0^i\|_{B^{\sigma-1}_{p,\infty}}\les \|u_0\|^2_{B^{\sigma}_{p,\infty}}.
\end{align*}
Thus, picking $t2^{kn}\approx\ep$ with small $\ep$, we have
\bbal
\|\mathbf{S}_{t}(u_0)-u_0\|_{B^\sigma_{p,\infty}}\geq c\ep-C\ep^2\geq \tilde{c}\ep.
\end{align*}
This completes the proof of Theorem \ref{th1}.

\section*{Acknowledgements}
J. Li is supported by the National Natural Science Foundation of China (11801090 and 12161004) and Jiangxi Provincial Natural Science Foundation (20212BAB211004). Y. Yu is supported by the National Natural Science Foundation of China (12101011) and Natural Science Foundation of Anhui Province (1908085QA05). W. Zhu is supported by the Guangdong Basic and Applied Basic Research Foundation (2021A1515111018).

\section*{Data Availability} Data sharing is not applicable to this article as no new data were created or analyzed in this study.

\section*{Conflict of interest}
The authors declare that they have no conflict of interest.

\end{document}